# ON ADIC GENUS, POSTNIKOV CONJUGATES, AND LAMBDA-RINGS

DONALD YAU

ABSTRACT. Sufficient conditions on a space are given which guarantee that the $K$-theory ring and the ordinary cohomology ring with coefficients over a principal ideal domain are invariants of, respectively, the adic genus and the SNT set. An independent proof of Notbohm's theorem on the classification of the adic genus of $BS^3$ by $KO$-theory $\lambda$-rings is given. An immediate consequence of these results about adic genus is that the power series ring $\mathbf{Z}[[x]]$ admits uncountably many pairwise non-isomorphic $\lambda$-ring structures.

## 1. INTRODUCTION AND STATEMENT OF RESULTS

In this paper we study some problems about the adic genus and the SNT set of a space. Let us first recall the relevant definitions. For a nilpotent finite type space $X$, we denote by, respectively, $\mathrm{SNT}(X)$ and $\widehat{G}_0(X)$ the set of homotopy types of spaces with the same $n$-type as $X$ for all $n$ and the set of homotopy types of nilpotent finite type spaces $Y$ with $\widehat{Y} \simeq \widehat{X}$ and $Y_0 \simeq X_0$. Here $\widehat{Y}$ denotes $\prod_p Y_p^\wedge$ with $Y_p^\wedge$ the $p$-completion of $Y$ and $Y_0$ denotes the rationalization of $Y$ (see [6]), and two spaces $X$ and $Y$ are said to have the same $n$-type if their Postnikov approximations through dimension $n$ are homotopy equivalent. The sets $\widehat{G}_0(X)$ and $\mathrm{SNT}(X)$ are called the *adic genus* and the SNT *set* of $X$, respectively. Two spaces with the same $n$-type for all $n$ are also known as *Postnikov conjugates* [8]. Note that by Wilkerson [25] $\mathrm{SNT}(X) \subseteq \widehat{G}_0(X)$, at least when $X$ is connected nilpotent of finite type.

The adic genus of an infinite dimensional space $X$ is often a very big set. For instance, Møller [17] proved that whenever $G$ is a compact connected non-abelian Lie group, the (adic) genus of its classifying space $BG$ is uncountably large. It is an important problem to find computable homotopy invariants which can distinguish between spaces with the same adic genus. Such a result was achieved by Notbohm and Smith [19] (see also [11]). Recall that $X$ is said to be an *adic fake Lie group of type $G$* if (1) $X = \Omega BX$ is a finite loop space and (2) $BX$ lies in the adic genus of $BG$. Notbohm and Smith showed that if $X$ is an adic fake Lie group of type $G$, where $G$ is a simply-connected compact Lie group, then $BX$ is homotopy equivalent to $BG$ if and only if $K(BX) \cong K(BG)$ as $\lambda$-rings. Here $K(Y)$ denotes the complex $K$-theory of $Y$.

With this result in mind, a natural question is then the following:

> *Is it really necessary to take into account the $\lambda$-operations in order to distinguish $BG$ from its fakes?*

Our first main result shows that for a large class of spaces this question has a positive answer.

**Theorem 1.1.** *Let $X$ be a simply-connected space of finite type whose integral homology is torsion free and is concentrated in even dimensions and whose $K$-theory filtered ring is a finitely generated power series ring over $\mathbf{Z}$. If $Y$ belongs to the adic genus of $X$, then there exists a filtered ring isomorphism $K(X) \cong K(Y)$.*

A couple of remarks are in order. First, this theorem has a variant in which complex $K$-theory (resp. $\mathbf{Z}$) is replaced with $KO^*$-theory (resp. $KO^* = KO^*(\text{pt})$), provided the integral homology of $X$ (which is simply-connected of finite type) is torsion free and is concentrated in dimensions divisible by 4 (e.g. $X = B\text{Sp}(n)$). This variant admits a proof which is essentially identical with that of Theorem 1.1.

Second, Theorem 1.1 and the result of Notbohm and Smith alluded to above imply that within the adic genus of $X = BG$ which satisfies the hypotheses of Theorem 1.1, there is only one $K$-theory ring (up to isomorphism) but there are at least two non-isomorphic $K$-theory $\lambda$-rings. This phenomenon leads us to the following questions.

> *How diverse is $K$-theory (or $KO$-theory) $\lambda$-ring in the adic genus of a space? How many $\lambda$-ring structures can a given power series ring support?*

Regarding the first question we found that in the case $X = BS^3$, the $KO$-theory $\lambda$-ring is as diverse as it can be.

**Theorem 1.2.** *Let $X$ and $Y$ be spaces in the adic genus of $BS^3$. Then $X$ and $Y$ are homotopy equivalent if and only if there exists a filtered $\lambda$-ring isomorphism $KO^*(X) \cong KO^*(Y)$.*

A more general result was obtained by Notbohm [18] by using topological realization properties of self homomorphisms of $K$-theory $\lambda$-rings of classifying spaces. Our proof of Theorem 1.2 is independent of Notbohm's and uses directly Rector's classification of the genus of $BS^3$ [20] and the $KO$-analogue of Theorem 1.1 when $X = BS^3$.

Before stating our next result, which concerns $\lambda$-ring structures over a given ring, let us first recall some relevant definitions. The reader is referred to [3, 13] for more information about $\lambda$-rings. A $\lambda$-*ring* is a commutative ring $R$ with unit together with functions $\lambda^i \colon R \to R$ ($i = 0, 1, \ldots$) satisfying certain properties. Given a $\lambda$-ring $(R, \{\lambda^i\})$ the *Adams operations* $\psi^k \colon R \to R$ ($k = 1, 2, \ldots$) are defined inductively by the Newton formulae:
$$\psi^k(a) - \lambda^1(a)\psi^{k-1}(a) + \cdots + (-1)^{k-1}\lambda^{k-1}(a)\psi^1(a) = (-1)^{k-1}k\lambda^k(a).$$
The Adams operations satisfy the following properties.
  1. All the $\psi^k \colon R \to R$ are ring homomorphisms.
  2. $\psi^1 = \text{Id}$ and $\psi^k \psi^l = \psi^{kl}$ for all $k, l \geq 1$.
  3. $\psi^p(a) \equiv a^p \mod pR$ for all primes $p$ and all $a$ in $R$.

On the other hand, if $S$ is a *torsion free* commutative ring with unit together with maps $\psi^k \colon S \to S$ ($k = 1, 2, \ldots$) satisfying the above three properties, then a theorem of Wilkerson [27] says that there exists a unique $\lambda$-ring structure over $S$ with the given $\psi^k$ as Adams operations.

Now since the adic genus of $BS^3$ is uncountable [20], Theorems 1.1 and 1.2 imply that, even in the one-variable case, a power series ring can support uncountably many pairwise non-isomorphic $\lambda$-ring structures.

**Corollary 1.3.** *There exist uncountably many distinct isomorphism classes of $\lambda$-ring structures over the power series ring $\mathbf{Z}[[x]]$.*

The analogous question of how many $\lambda$-ring structures the *polynomial ring* $\mathbf{Z}[x]$ supports has been studied by Clauwens [7]. Employing the theory of commuting polynomials, he showed that there are essentially only two non-isomorphic $\lambda$-ring structures on the polynomial ring $\mathbf{Z}[x]$. We, however, have not been able to establish any connections between Clauwens' result and our Corollary 1.3.

The last main result of this note concerns the SNT set of a space. There are many interesting spaces whose SNT sets are nontrivial (in fact, uncountable); see, for example, [9, 15, 22, 23]. Several homotopy invariants which can distinguish between spaces of the same $n$-type for all $n$ have been found. These include [16] $\text{Aut}(-)$ and $\text{WI}(-)$, the group of homotopy classes of homotopy self-equivalences and the group of weak identities, and [24] $\text{End}(-)$, the monoid of homotopy classes of self-maps. One might wonder if there are more familiar and more computable homotopy invariants, such as integral cohomology and $K$-theory, which can distinguish between spaces in the same SNT set. Of course, our Theorem 1.1 together with [25] imply that if $X$ is as in Theorem 1.1 and if $Y$ lies in $\mathrm{SNT}(X)$, then the integral cohomology (resp. $K$-theory) rings of $X$ and $Y$ are isomorphic.

   *Does this hold for more general spaces?*

The following result gives an answer to this question for the ordinary cohomology case.

**Theorem 1.4.** *Let $\Lambda$ be a principal ideal domain. Let $X$ be a connected space of finite type whose ordinary cohomology ring with coefficients over $\Lambda$ is a finitely generated graded algebra over $\Lambda$. If $Y$ has the same n-type as $X$ for all n, then there exists an isomorphism $H^*(X; \Lambda) \cong H^*(Y; \Lambda)$ of cohomology rings.*

It should be emphasized that in Theorem 1.4 the hypothesis $Y \in \mathrm{SNT}(X)$ *cannot* be weakened to the condition $Y \in \widehat{G}_0(X)$. Indeed, in [5] Bokor constructed two spaces, both two cell complexes, with the same genus but whose integral cohomology rings are non-isomorphic.

This finishes the presentation of the results in this paper. The rest of the paper is organized as follows. Section 2 contains preliminary materials. The proofs of Theorems 1.1, 1.2, and 1.4 are given, one in each section, in §3 - §5.

## Acknowledgement

The results in this paper constitute part of the author's doctoral research in progress. The author would like to express his sincerest gratitude to his advisor Professor Haynes Miller for superb guidance and encouragement.

## 2. Preliminaries

A result that we need is a straightforward generalization of Wilkerson's Classification Theorem [25, Thm. 1.1] of spaces of the same $n$-type for all $n$, as given in [28, Thm. 12]. We will use it in two algebraic categories: the category of commutative filtered rings (for Theorem 1.1) and the category of complete graded algebras over $\Lambda$ (a principal ideal domain) (for Theorem 1.4). We will now recall the relevant notations and terminology.

Let $\mathcal{M}$ be a pointed category. A *localization functor* on $\mathcal{M}$ consists of:

1. A functor $L\colon \mathcal{M} \to \mathcal{M}$.
2. A natural transformation $\eta\colon \mathrm{Id} \to L$ (where Id denotes the identity functor on $\mathcal{M}$).
3. A natural equivalence $L^2 \cong L$.

We will abuse notation and refer to $L$ as the localization functor. An object $X$ in $\mathcal{M}$ is called *L-local* if the coaugmentation $\eta_X\colon X \to LX$ is an isomorphism. A *localization tower* $\{L_n\}_{n\geq 0} = \{L_0 \leftarrow L_1 \leftarrow \cdots\}$ on $\mathcal{M}$ consists of localization functors $L_n$ ($n \geq 0$) such that every $L_n$-local object is also $L_{n+1}$-local and the natural transformations $L_{n+1} \to L_n$ ($n \geq 0$) are induced by the universal property of $L_{n+1}$-localization. The reader is referred to Adams' book [1] for more information on localization functors.

Now suppose that $T = \{L_n\}_{n\geq 0}$ is a localization tower on $\mathcal{M}$ and that $X$ is an object in $\mathcal{M}$. Following Wilkerson [25] we denote by $\mathrm{SNT}(X)$ the collection of isomorphism classes of objects $Y$ in $\mathcal{M}$ with $L_nX \cong L_nY$ for all $n \geq 0$, and by $\mathrm{SNT}(X)^\wedge$ the subset of $\mathrm{SNT}(X)$ consisting of those (isomorphism classes of) objects $Y$ such that the natural map $Y \to \varprojlim L_nY$ is an isomorphism. An object $Y$ which satisfies the last condition is said to be *complete*. Note that $\mathrm{SNT}(X)$ is a pointed set with $[X]$ as its distinguished element. For an object $X$ in $\mathcal{M}$ its automorphism group is denoted by $\mathrm{Aut}(X)$. Finally, if $\{G_n\}_{n\geq 0}$ is an inverse system of (not-necessarily abelian) groups, then the pointed set $\varprojlim{}^1 G_n$ is defined as [6] the orbit set of the action

$$\prod_n G_n \times \prod_n G_n \to \prod_n G_n$$

$$((\alpha_n), (\beta_n)) \longmapsto (\alpha_n \beta_n \rho_{n+1}(\alpha_{n+1})^{-1})$$

where $\rho_{n+1}\colon G_{n+1} \to G_n$ is the structure map.

As noted in [28] Wilkerson's proof of [25, Thm. 1.1] can be easily adapted to yield the following classification result.

**Theorem 2.1** (Wilkerson [25]). *Let $X$ be an object in $\mathcal{M}$ which is complete with respect to a localization tower $T = \{L_n\}_{n\geq 0}$. Then there exists a bijection $\varprojlim{}^1 \mathrm{Aut}(L_nX) \cong \mathrm{SNT}(X)^\wedge$ of pointed sets.*

In particular, if the tower $\{\mathrm{Aut}(L_nX)\}$ is surjective, i.e. the maps $\mathrm{Aut}(L_{n+1}X) \to \mathrm{Aut}(L_nX)$ are surjective for all $n$ sufficiently large, then $\mathrm{SNT}(X)^\wedge = \{[X]\}$.

3. Proof of Theorem 1.1

In this section we prove Theorem 1.1.

We will work within the category **CFR** of commutative filtered rings. Thus, an object in **CFR** is a pair $(R, \{I^n\}_{n \geq 0})$ consisting of a commutative ring $R$ with unit and a decreasing filtration $R = I^0 \supseteq I^1 \supseteq \cdots$ of ideals. Whenever there is no danger of confusion we will not mention the filtration explicitly. Every space $Z$ (with the homotopy type of a CW complex) gives rise to an object $(K(Z), \{K_n(Z)\}_{n \geq 0})$ in **CFR**, where $K(Z)$ denotes the complex $K$-theory ring of $Z$ and $K_n(Z) = \ker(i_n^* \colon K(Z) \to K(Z_{n-1}))$ with $i_n \colon Z_{n-1} \subset Z$ the inclusion of the $(n-1)$st skeleton.

We will apply Theorem 2.1 with $\mathcal{M} = $ **CFR** and with the following localization tower. For $j \geq 0$ define the functor

$$L_j \colon \mathbf{CFR} \to \mathbf{CFR}, \quad L_j R = R/I^j$$

where $R/I^j$ has the induced quotient filtration. Each functor $L_j$ comes equipped with a natural coaugmentation $\mathrm{Id} \to L_j$, and it is easy to check that $\{L_j\}_{j \geq 0}$ is a localization tower on the category **CFR** in the sense of §2. Thus, given a commutative filtered ring $(R, \{I^j\})$ the notation $\mathrm{SNT}(R)$ now denotes the set of isomorphism classes of commutative filtered rings $(S, \{J^j\})$ with $R/I^j \cong S/J^j$ in **CFR** for all $j \geq 0$.

From this point on in this section, $X$ and $Y$ are as in Theorem 1.1. To prove Theorem 1.1 we need the following lemmas.

**Lemma 3.1.** *The natural map $K(Y) \to \varprojlim K(Y)/K_j(Y)$ is an isomorphism.*

**Lemma 3.2.** *The filtered ring $K(Y)$ lies in $\mathrm{SNT}(K(X))$.*

By hypothesis the $K$-theory filtered ring of $X$ has the form

$$K(X) = \mathbf{Z}[[c_1, \ldots, c_n]]$$

in which the generators $c_i$ are algebraically independent over $\mathbf{Z}$ and $K_j(X)$ is the ideal generated by the monomials of degrees at least $j$. Suppose that $c_i \in K_{d_i}(X)$, and let $N$ be an integer that is strictly greater than $d_i$ for $i = 1, \ldots, n$.

**Lemma 3.3.** *Every map in the inverse system $\{\mathrm{Aut}(L_j K(X))\}_{j > N}$ is surjective.*

Assuming these three lemmas for the moment let us give the proof of Theorem 1.1.

*Proof of Theorem 1.1.* Lemmas 3.1 and 3.2 combine to imply that $K(Y)$ lies in $\mathrm{SNT}(K(X))^{\wedge}$, which by Theorem 2.1 and Lemma 3.3 is the one-point set $\{[K(X)]\}$. Thus, $K(X) \cong K(Y)$ in **CFR**, as desired. □

Now we give the proofs of the lemmas above.

*Proof of Lemma 3.1.* According to [4, 2.5 and 7.1] the natural map $K(Y) \to \varprojlim K(Y)/K_j(Y)$ is an isomorphism if

$$(3.1) \qquad \varprojlim{}^1 E_r^{p,q} = 0 \quad \text{for all } (p,q),$$

where $E_r^{*,*}$ is the $E_r$-term in the $K^*$-Atiyah-Hirzebruch spectral sequence (AHSS) for $Y$ and $\varprojlim{}^1$ denotes the first right derived functor of the inverse limit functor. This is the case when the AHSS degenerates at the $E_2$-term. Thus, to prove (3.1) it suffices to show that $H^*(Y;\mathbf{Z})$ is concentrated in even dimensions because in that case there is no room for differentials in the AHSS. So pick an odd integer $N$. We must show that

$$(3.2) \qquad H^N(Y;\mathbf{Z}) \cong 0.$$

By the Universal Coefficient Theorem it suffices to show that the integral homology of $Y$ is torsion free and is concentrated in even dimensions, which hold because $Y \in \widehat{G}_0(X)$.

This finishes the proof of Lemma 3.1. □

*Proof of Lemma 3.2.* We have to show that for each $j > 0$ there is an isomorphism

$$(3.3) \qquad K(Y)/K_j(Y) \cong K(X)/K_j(X)$$

in **CFR**. It follows from the hypothesis that for each $j > 0$, $K(Y)/K_j(Y) \in \widehat{G}_0(K(X)/K_j(X))$ where $\widehat{G}_0(R)$ for $R \in$ **CFR** is defined in terms of $R \otimes \mathbf{Q}$ and $R \otimes \widehat{\mathbf{Z}}$ in exactly the same way the adic genus of a space is defined. Now Wilkerson's proof of [26, Thm. 3.8] can be easily adapted to show that there is a surjective map $\mathrm{Caut}((K(X)/K_j(X)) \otimes \mathbf{Q} \otimes \widehat{\mathbf{Z}}) \to \widehat{G}_0(K(X)/K_j(X))$ (see [26, §3] for the definition of $\mathrm{Caut}(-)$). Then [26, 3.7] implies that the image of this map is constant because $K(X)/K_j(X)$ is a finitely generated abelian group. This implies immediately (3.3).

This finishes the proof of Lemma 3.2. □

*Proof of Lemma 3.3.* Fix an integer $j > N$ and pick a filtered ring automorphism $\sigma$ of $K(X)/K_j(X)$. We must show that $\sigma$ can be lifted to a filtered ring automorphism of $K(X)/K_{j+1}(X)$. For $1 \leq i \leq n$ pick any lift of $\sigma(c_i)$ to $K(X)/K_{j+1}(X)$ and define $\hat\sigma(c_i) \in K(X)/K_{j+1}(X)$ to be such a lift. Since there are no relations among the $c_i$ in $K(X)$ it is easy to see that $\hat\sigma$ gives rise to a well-defined filtered ring endomorphism of $K(X)/K_{j+1}(X)$, and it will be a desired lift of $\sigma$ once it is shown to be bijective.

To show that $\hat\sigma \colon K(X)/K_{j+1}(X) \to K(X)/K_{j+1}(X)$ is surjective, it suffices to show that each $c_i \in K(X)/K_{j+1}(X)$ lies in the image of $\hat\sigma$, since $K(X)/K_{j+1}(X)$ is generated as a filtered ring by the $c_i$. So fix an integer $i$ with $1 \leq i \leq n$. We know that there exists an element $g_i \in K(X)/K_j(X)$ such that

$$(3.4) \qquad \sigma(g_i) = c_i.$$

Pick any lift of $g_i$ to $K(X)/K_{j+1}(X)$, call it $g_i$ again, and observe that (3.4) implies that

$$(3.5) \qquad \hat\sigma(g_i) = c_i + \alpha_i$$

in $K(X)/K_{j+1}(X)$ for some $\alpha_i \in K_j(X)/K_{j+1}(X)$. We will alter $g_i$ to obtain a $\hat\sigma$-pre-image of $c_i$ as follows. The ideal $K_j(X)/K_{j+1}(X)$ is generated by certain monomials in $c_1, \ldots, c_n$, namely,

$$c_1^{i_1} \cdots c_n^{i_n}, \quad (i_1, \ldots, i_n) \in J_j$$

where $J_j$ is the set of ordered $n$-tuples $(i_1, \ldots, i_n)$ of nonnegative integers with $\sum_{l=1}^n d_l i_l = j$. Thus, for every $(i_1, \ldots, i_n) \in J_j$ there exists an integer $a_{(i_1, \ldots, i_n)}$ such that

$$\alpha_i = \sum_{(i_1, \ldots, i_n) \in J_j} a_{(i_1, \ldots, i_n)} \, c_1^{i_1} \cdots c_n^{i_n}.$$

Now define $\bar g_i$ in $K(X)/K_{j+1}(X)$ by

$$\bar g_i \equiv g_i - \sum_{(i_1, \ldots, i_n) \in J_j} a_{(i_1, \ldots, i_n)} \, g_1^{i_1} \cdots g_n^{i_n}.$$

We claim that

$$(3.6) \qquad \hat\sigma(\bar g_i) = c_i \in K(X)/K_{j+1}(X).$$

It clearly suffices to prove that

$$(3.7) \qquad \hat\sigma(g_1^{i_1} \cdots g_n^{i_n}) = c_1^{i_1} \cdots c_n^{i_n} \in K(X)/K_{j+1}(X)$$

for each $(i_1, \ldots, i_n) \in J_j$. Now in $K(X)/K_{j+1}(X)$ one has

$$\begin{aligned}
\hat{\sigma}(g_1^{i_1} \cdots g_n^{i_n}) &= \hat{\sigma}(g_1)^{i_1} \cdots \hat{\sigma}(g_n)^{i_n} \\
&= (c_1 + \alpha_1)^{i_1} \cdots (c_n + \alpha_n)^{i_n} \quad \text{by (3.5)} \\
&= c_1^{i_1} \cdots c_n^{i_n} + (\text{terms of filtration } > j) \\
&= c_1^{i_1} \cdots c_n^{i_n}.
\end{aligned}$$

This proves (3.7), and hence (3.6), and therefore $\hat{\sigma}$ is surjective.

It remains to show that $\hat{\sigma}$ is injective. Since any surjective endomorphism of a finitely generated abelian group is also injective and since $K(X)/K_{j+1}(X)$ is a finitely generated abelian group, it follows that $\hat{\sigma}$ is injective as well. Thus, $\hat{\sigma}$ is an automorphism of $K(X)/K_{j+1}(X)$ and is a lift of $\sigma$.

This finishes the proof of Lemma 3.3. □

## 4. Proof of Theorem 1.2

In this section we prove Theorem 1.2. The arguments in this section, especially Lemma 4.2 below, are inspired by Rector's [21, §4].

We begin by noting that an argument entirely similar to the proof of Theorem 1.1 implies that whenever $X$ belongs to $\widehat{G}_0(BS^3)$,

$$KO^*(X) \cong KO^*[[x]]$$

as filtered rings, where $x \in KO_4^4(X)$ is a representative of an integral generator $x_4 \in H^4(X; \mathbf{Z}) = E_2^{4,0}$ in the $KO^*$-Atiyah-Hirzebruch spectral sequence for $X$. Here $KO_b^a(X)$ denotes the subgroup of $KO^a(X)$ consisting of elements $u$ which restrict to 0 under the natural map $KO^a(X) \to KO^a(X_{b-1})$. Such an element $u$ is said to be in *degree a* and *filtration b*.

Now we recall the relevant notations, definitions, and results regarding Rector's classification of the genus of $BS^3$ [20]. Let $\xi \in \pi_{-4}KO$ and $b_R \in \pi_{-8}KO$ be the generators so that $\xi^2 = 4b_R$. As usual, denote by $\psi^k$ ($k = 1, 2, \ldots$) the Adams operations. Since $\Omega X \simeq S^3$ it follows as in [21, §4] that there exists an integer $a$, depending on the choice of the representative $x$, such that:

1. $\psi^2(\xi x) = 4\xi x + 2ab_R x^2 \pmod{KO_9^0(X)}$.
2. The integer $a$ is well-defined (mod 24). This means that if $x'$ is another representative of $x_4$ with corresponding integer $a'$, then $a \equiv a' \pmod{24}$, and if $x_4$ is replaced with $-x_4$, then $a$ will be replaced with $-a$. We can therefore write $a(X)$ for $a$.
3. $a(X) \equiv \pm 1, \pm 5, \pm 7,$ or $\pm 11 \pmod{24}$.

The last condition above follows from the examples constructed by Rector in [21, §5] and James' result [10] which says that there are precisely eight homotopy classes of homotopy-associative multiplications on $S^3$. These eight classes can be divided into four pairs with each pair consisting of a homotopy class of multiplication and its inverse.

Rector's invariant $(X/p)$ for $p$ an odd prime is defined as follows [20]. The Adem relation $P^1 P^1 = 2P^2$ implies that

$$P^1 \overline{x}_4 = \pm 2\overline{x}_4^{(p+1)/2}$$

in $H^*(X; \mathbf{Z}/p)$, where $\overline{x}_4$ is the mod $p$ reduction of the integral generator $x_4$. Then $(X/p) \in \{\pm 1\}$ is defined as the sign on the right-hand side of this equation.

The invariant $(X/2)$ and a canonical choice of orientation of the integral generator are given as follows. Using the mod 24 integer $a(X)$, define

$$(4.1) \qquad ((X/2), (X/3)) = \begin{cases} (1, 1) & \text{if } a(X) \equiv \pm 1 \bmod 24; \\ (1, -1) & \text{if } a(X) \equiv \pm 5 \bmod 24; \\ (-1, 1) & \text{if } a(X) \equiv \pm 7 \bmod 24; \\ (-1, -1) & \text{if } a(X) \equiv \pm 11 \bmod 24. \end{cases}$$

The orientation is then chosen so that $(X/3)$ is as given in (4.1). This definition of Rector's invariants coincides with the original one (cf. [14, §9]). Now we can recall the classification theorem of the genus of $BS^3$ [20].

**Theorem 4.1** (Rector). *The $(X/p)$ for $p$ primes provide a complete list of classification invariants for the (adic) genus of $BS^3$. Any combination of values of the $(X/p)$ can occur. If $X = \mathbf{HP}^\infty$ then $(X/p) = 1$ for all primes $p$.*

From now on in this section, $X$ and $Y$ are as in Theorem 1.2. Now we can prove Theorem 1.2.

*Proof of Theorem 1.2.* In view of Theorem 4.1, to prove Theorem 1.2 it suffices to show that if there exists a filtered $\lambda$-ring isomorphism $\sigma \colon KO^*(X) \cong KO^*(Y)$, then $a(X) \equiv \pm a(Y) \pmod{24}$ and $(X/p) = (Y/p)$ for all odd primes $p$. We will prove these in Lemmas 4.2 and 4.3 below, thereby proving Theorem 1.2. □

As explained above, $KO^*(X) = KO^*[[x]]$ and $KO^*(Y) = KO^*[[y]]$ with $x \in KO_4^4(X)$ and $y \in KO_4^4(Y)$ representing, respectively, the integral generators $x_4 \in H^4(X; \mathbf{Z})$ and $y_4 \in H^4(Y; \mathbf{Z})$.

**Lemma 4.2.** *If there exists a filtered $\lambda$-ring isomorphism $\sigma \colon KO^*(X) \cong KO^*(Y)$, then $a(X) \equiv \pm a(Y) \pmod{24}$.*

*Proof.* Since $\sigma$ is a ring isomorphism, we have
$$\sigma(\xi x) = \varepsilon \xi y + \sigma_2 b_R y^2 \pmod{KO_9^0(Y)}$$
for some integer $\sigma_2$ and $\varepsilon \in \{\pm 1\}$. Computing modulo $KO_9^0(X)$ we have
$$4\sigma(b_R x^2) = \sigma(\xi x)^2 = \xi^2 y^2 = 4 b_R y^2,$$
and therefore $\sigma(b_R x^2) = b_R y^2 \pmod{KO_9^0(Y)}$. First we claim that there is an equality
$$(4.2) \qquad a(X) = 6\sigma_2 + \varepsilon a(Y).$$

To prove (4.2) we will compute both sides of $\sigma \psi^2(\xi x) = \psi^2 \sigma(\xi x)$ (mod $KO_9^0(Y)$). Working modulo $KO_9^0(Y)$ we have, on the one hand,
$$\begin{aligned}
\sigma \psi^2(\xi x) &= \sigma(4\xi x + 2a(X) b_R x^2) \\
&= 4(\varepsilon \xi y + \sigma_2 b_R y^2) + 2a(X) b_R y^2 \\
&= 4\varepsilon \xi y + (4\sigma_2 + 2a(X)) b_R y^2.
\end{aligned}$$
On the other hand, still working modulo $KO_9^0(Y)$, we have
$$\begin{aligned}
\psi^2 \sigma(\xi x) &= \varepsilon \psi^2(\xi y) + \sigma_2 \psi^2(b_R y^2) \\
&= \varepsilon(4\xi y + 2a(Y) b_R y^2) + \sigma_2(2^4 b_R y^2) \\
&= 4\varepsilon \xi y + (16\sigma_2 + 2\varepsilon a(Y)) b_R y^2.
\end{aligned}$$
Equation (4.2) now follows by equating the coefficients of $b_R y^2$.

In view of (4.2), to finish the proof of Lemma 4.2 it is enough to establish
$$(4.3) \qquad \sigma_2 \equiv 0 \pmod{4}.$$
To prove (4.3), note that since $\sigma$ is a $KO^*$-module map, we have $\xi \sigma(x) = \sigma(\xi x)$. Since $\sigma$ is a ring isomorphism, we also have
$$\sigma(x) = \varepsilon' y + \sigma_2' \xi y^2 \pmod{KO_9^4(Y)}$$
for some integer $\sigma_2'$ and $\varepsilon' \in \{\pm 1\}$. Therefore, working modulo $KO_9^0(Y)$ we have
$$\begin{aligned}
\xi \sigma(x) &= \varepsilon' \xi y + \sigma_2' \xi^2 y^2 \\
&= \varepsilon' \xi y + 4 \sigma_2' b_R y^2 \\
&= \varepsilon \xi y + \sigma_2 b_R y^2.
\end{aligned}$$
In particular, by equating the coefficients of $b_R y^2$ we obtain $\sigma_2 = 4\sigma_2'$, thereby proving (4.3).

This completes the proof of Lemma 4.2. $\square$

**Lemma 4.3.** *If there exists a filtered $\lambda$-ring isomorphism $\sigma \colon KO^*(X) \cong KO^*(Y)$, then $(X/p) = (Y/p)$ for each odd prime $p$.*

*Proof.* It follows from Theorem 1.1 that $K^*(X) \cong K^*[[u_x]]$ with $u_x \in K_4^4(X)$ a representative of the integral generator $x_4 \in H^4(X; \mathbf{Z}) = E_2^{4,0}$ in the $K^*$-Atiyah-Hirzebruch spectral sequence. Moreover, we may choose $u_x$ so that $c(x) = u_x$, where $c \colon KO^*(X) \to K^*(X)$ is the complexification map. Similar remarks apply to $Y$ so that $K^*(Y) \cong K^*[[u_y]]$.

Now denote by $b \in \pi_{-2}K$ the Bott element and let $p$ be a fixed odd prime. We first claim that

$$(4.4) \quad \psi^p(b^2 u_x) = (b^2 u_x)^p + 2(X/p)\, p\, (b^2 u_x)^{(p+1)/2} + p\, w_x + p^2 x_0$$

for some $w_x \in K_{2p+3}^0(X)$ and some $x_0 \in K_4^0(X)$. To see this, note that since $b^2 u_x \in K_4^0(X)$, it follows from Atiyah's theorem [2, 5.6] that

$$\psi^p(b^2 u_x) = (b^2 u_x)^p + p\, x_1 + p^2 x_0$$

for some $x_i \in K_{4+2i(p-1)}^0(X)$ ($i = 0, 1$); moreover, $\overline{x}_1 = P^1 \overline{b^2 u_x}$, where $\overline{z}$ is the mod $p$ reduction of $z$. Thus, to prove (4.4) it is enough to show that

$$(4.5) \quad x_1 = 2(X/p)(b^2 u_x)^{(p+1)/2} + w_x + p\, z_x$$

for some $w_x \in K_{2p+3}^0(X)$ and some $z_x \in K_{2p+2}^0(X)$. Now in $H^*(X; \mathbf{Z}) \otimes \mathbf{Z}/p$ we have

$$\overline{x}_1 = P^1 \overline{b^2 u_x} = P^1 \overline{x}_4 = 2(X/p)\, \overline{x}_4^{(p+1)/2} = 2(X/p)\, \overline{b^2 u_x}^{(p+1)/2},$$

from which (4.5) follows immediately. As remarked above, this also establishes (4.4).

Now the $\lambda$-ring isomorphism $\sigma$ induces via $c$ a $\lambda$-ring isomorphism $\sigma_c \colon K^*(X) \cong K^*(Y)$. By composing $\sigma_c$ with a suitable $\lambda$-ring automorphism of $K^*(Y)$ if necessary, we obtain a $\lambda$-ring isomorphism $\alpha \colon K^*(X) \cong K^*(Y)$ with the property that

$$(4.6) \quad \alpha(b^2 u_x) = b^2 u_y + \text{ higher terms in } b^2 u_y.$$

Using (4.4) and (4.6) it is then easy to check that

$$(4.7) \quad \alpha \psi^p(b^2 u_x) = 2(X/p)\, p\, (b^2 u_y)^{(p+1)/2} \quad (\text{mod } K_{2p+3}^0(Y) \text{ and } p^2)$$

and

$$(4.8) \quad \psi^p \alpha(b^2 u_x) = 2(Y/p)\, p\, (b^2 u_y)^{(p+1)/2} \quad (\text{mod } K_{2p+3}^0(Y) \text{ and } p^2).$$

Since $\alpha \psi^p = \psi^p \alpha$ it follows from (4.7) and (4.8) that

$$2(X/p)\, p \equiv 2(Y/p)\, p \pmod{p^2},$$

or, equivalently,
$$2(X/p) \equiv 2(Y/p) \pmod{p}.$$
But $p$ is assumed odd, and so $(X/p) \equiv (Y/p) \pmod{p}$. Hence $(X/p) = (Y/p)$, as desired.

This finishes the proof of Lemma 4.3. □

## 5. Proof of Theorem 1.4

In this final section we give the proof of Theorem 1.4, which is structurally similar to the proof of Theorem 1.1. All cohomology groups and rings have coefficients over a fixed principal ideal domain $\Lambda$. Consider the category **CGA** of complete graded algebras over $\Lambda$. Thus, an object in **CGA** is a $\Lambda$-algebra $R = \prod_{0 \leq i} R^i$ with $R^i$ the homogeneous part of degree $i$. Note that the natural map $R \to \varprojlim R^{\leq n}$ is an isomorphism. Every space $Z$ gives rise naturally to an object $H^{**}(Z) = \prod_{0 \leq i} H^i(Z)$ in **CGA**.

We will use Theorem 2.1 with $\mathcal{M} = \mathbf{CGA}$ and with the following localization tower. For every nonnegative integer $j$ define the functor
$$L_j \colon \mathbf{CGA} \to \mathbf{CGA} \quad \text{by} \quad L_j R = R^{\leq j} = \frac{R}{\prod_{j<i} R^i}.$$
Each functor $L_j$ comes equipped with a natural coaugmentation $\mathrm{Id} \to L_j$ and it is easy to check that $\{L_j\}_{j \geq 0}$ is a localization tower on **CGA** in the sense of §2. Given an object $R$ in **CGA** the notation $\mathrm{SNT}(R)$ now denotes the set of isomorphism classes of complete graded algebras $S$ over $\Lambda$ such that $L_j R \cong L_j S$ in **CGA** for all $j \geq 0$. (Note that every object in **CGA** is complete with respect to this localization tower, and so $\mathrm{SNT}(-) = \mathrm{SNT}(-)^\wedge$.)

From now on $X$ and $Y$ are as in Theorem 1.4. Observe that since the cohomology ring $H^*(Z)$ is the associated graded of $H^{**}(Z)$ with respect to the filtration $\{H^{\leq j}(X)\}$, Theorem 1.4 is an immediate consequence of the following result.

**Proposition 5.1.** *There exists an isomorphism $H^{**}(X) \cong H^{**}(Y)$ in* **CGA**.

To prove Proposition 5.1 we need the following lemmas.

**Lemma 5.2.** *The complete graded $\Lambda$-algebra $H^{**}(Y)$ belongs to* $\mathrm{SNT}(H^{**}(X))$.

By hypothesis the cohomology ring of $X$ has the form
$$H^*(X) = \Lambda[x_1, \ldots, x_s]/J$$
for some homogeneous generators $x_i$ ($1 \leq i \leq s$) and some homogeneous ideal $J$. Suppose that $|x_i| = d_i$. Note that since the polynomial ring $\Lambda[x_1, \ldots, x_s]$ is Noetherian, the ideal $J$ is generated by finitely many polynomials, say, $f_m$ ($1 \leq m \leq k$). We can thus choose an integer $N$ which is strictly greater than the $d_i$ ($1 \leq i \leq s$) and the degree of any nontrivial monomial in any $f_m$ ($1 \leq m \leq k$).

**Lemma 5.3.** *Every map in the inverse system $\{\mathrm{Aut}\,(L_j\, H^{**}(X))\}_{j>N}$ is surjective.*

*Proof of Proposition 5.1.* Theorem 2.1 and Lemma 5.3 combine to imply that $\mathrm{SNT}\,(H^{**}(X)) = \{[H^{**}(X)]\}$. This together with Lemma 5.2 yield the desired isomorphism: $H^{**}(X) \cong H^{**}(Y)$. □

It remains to give the proofs of Lemmas 5.2 and 5.3.

*Proof of Lemma 5.2.* This follows immediately from the assumption that $X$ and $Y$ have the same $n$-type for all $n$. □

*Proof of Lemma 5.3.* Let $j$ be any integer strictly greater than $N$ and let $\sigma$ be an automorphism of the graded ring $H^{\leq j}(X)$. We want to show that $\sigma$ can be lifted to an automorphism of $H^{\leq j+1}(X)$. Now $\sigma(f_m)$ is 0 in $H^{\leq j}(X)$ for $m = 1, \ldots, k$. But since each $\sigma(x_i)$ is homogeneous of degree $d_i$, each monomial of $\sigma(f_m) = f_m(\sigma(x_1), \ldots, \sigma(x_s))$ still has degree less than $N$. Therefore, there is a well-defined graded ring map
$$\hat{\sigma} \colon H^{\leq j+1}(X) \to H^{\leq j+1}(X)$$
satisfying
$$\hat{\sigma}|H^{\leq j}(X) = \sigma \colon H^{\leq j}(X) \to H^{\leq j}(X).$$
We will be done once we show that $\hat{\sigma}$ is bijective.

It is clear that $\hat{\sigma}$ is surjective because $H^{\leq j+1}(X)$ is generated as an algebra by the $x_i$ ($1 \leq i \leq s$) and each $x_i$ is in the image of $\hat{\sigma}$, since this is true for $\sigma$.

Finally, since $H^{\leq j+1}(X)$ is a finitely generated $\Lambda$-module, the injectivity of $\hat{\sigma}$ now follows from its surjectivity. Therefore, $\hat{\sigma}$ is an automorphism of $H^{\leq j+1}(X)$ and is a lift of $\sigma$.

This finishes the proof of Lemma 5.3. □

Department of Mathematics, MIT Room 2-230, 77 Massachusetts Avenue, Cambridge, MA 02139, USA

*E-mail address*: `donald@math.mit.edu`